\documentclass[12pt,a4paper]{article}
\RequirePackage{amsmath,amssymb,amsfonts}
\usepackage[english,francais]{babel} 
\newcommand{\sth}{{{{\mbox{${\mbox{\scriptsize\boldmath$\theta$}}$}}}}}

\newcommand{\xbrt}{X^{{\rm br},\sth}}

\newcommand{\sfrac}[2]{{\textstyle\frac{#1}{#2}}}

\newcommand{\bp}{{\bf p}}

\newcommand{\bq}{{\bf q}}
\newcommand{\bw}{{\bf w}}

\newcommand{\bt}{{\bf t}}
\newcommand{\bJ}{{\bf J}}

\newcommand{\bT}{{\bf T}}

\newcommand{\btheta}{{\mbox{\boldmath$\theta$}}}
\newcommand{\bTheta}{{\mbox{\boldmath$\Theta$}}}
\newcommand{\Thfinite}{\bTheta_{{\rm finite}}}

\renewcommand{\sp}{\sigma({\bf p})}
\newcommand{\spn}{\sigma({\bf p}_n)}

\newcommand{\R}{\mathbb{R}}

\newcommand{\D}{\mathbb{D}}

\newcommand{\TT} { {\cal T }}

\newcommand{\cyc} { {\cal C }}

\newcommand{\com} { {\cal B }}
\def\build#1_#2^#3{\mathrel{
\mathop{\kern 0pt#1}\limits_{#2}^{#3}}}

\newcommand{\bb}{B^{|{\rm br}|}}
\newcommand{\TTth}{{\cal T}^{\sth}}

\newcommand{\HT}{{\rm ht}}

\def\cq{$\hfill \square$}
\def\un{\underline}
\def\d{{\rm d}}

\def\eps{\varepsilon}

\def\ba{\begin{eqnarray*}}
\def\ea{\end{eqnarray*}}

\def\wt{\widetilde}

\def\rem{\noindent{\bf Remark. }}

\def\proof{\noindent{\bf Proof. }}

\def\card{{\rm Card\,}}

\newtheorem{thm}{Theorem}
\newtheorem{lmm}{Lemma}
\newtheorem{prp}{Proposition}
\newtheorem{crl}{Corollary}

\title{Weak convergence  of random $\bp$-mappings and  the exploration process
of inhomogeneous continuum random trees}

\author{David Aldous\thanks{Department of Statistics, U.C. Berkeley CA 94720-3860, USA}, Gr\'egory Miermont\thanks{DMA, \'Ecole Normale
Sup\'erieure, 45 rue d'Ulm, F-75230 Paris Cedex 05 and Laboratoire de
Probabilit\'es et Mod\`eles Al\'eatoires, Universit\'e Paris VI. 
{\tt miermont@dma.ens.fr}}
\ and Jim Pitman\thanks{Department of Statistics, U.C. Berkeley CA 94720-3860, USA}}
\date{}

\begin{document}

\selectlanguage{english}

\maketitle

\begin{abstract}
We study the asymptotics of the $\bp$-mapping model of random mappings on
$[n]$  as $n$ gets  large, under  a large  class of  asymptotic regimes  for the
underlying distribution $\bp$. % including the cases when a finite number of 
We encode these random mappings in random walks which are shown to converge to
a functional of the exploration process of inhomogeneous random
trees, 
this exploration process being derived 
(Aldous-Miermont-Pitman 2003)
from a bridge with exchangeable increments.
Our setting generalizes previous results by allowing
a finite number of 
``attracting  points'' to emerge. %  thus  generalizing the  result  of a  previous
\end{abstract}



\newpage

\section{Introduction}

We study the asymptotic behavior as $n\to\infty$
of random elements of the set $[n]^{[n]}$ of {\em mappings} from $[n] = \{1,2,\ldots,n\}$ to
$[n]$. Given a probability measure $\bp=(p_1,\ldots,p_n)$ on
$[n]$, define a random mapping $M$ as follows: for each
$i\in[n]$, map $i$ to $j$ with probability $p_j$, independently over
different $i$'s, so that 
\begin{equation}\label{pmapping}
P(M=m)=\prod_{i\in[n]}p_{m(i)}\, ,\qquad m\in[n]^{[n]}. 
\end{equation}
The random mapping $M$ is called the
{\em $\bp$-mapping}. In what follows, we will not be concerned about keeping
track of the labels of the mapping's digraph, so we will suppose that the
probability $\bp$ is {\em ranked}, i.e.\ $p_1\geq p_2\geq\ldots\geq
p_n>0$.  %(the positivity assumption is to avoid heavy discussions later). 

Now consider a sequence of such probabilities $\bp_n=(p_{n1},\ldots,p_{nn})$. 
Weak convergence of the associated
%In \cite{almpimap02}, we have described the weak convergence of the associated
$\bp$-mappings $M_n$
as $n\to\infty$ has been studied when $\bp_n$ satisfies an asymptotic negligibility condition,
namely, letting $\spn=(\sum_{1\leq i\leq n}p_{ni}^2)^{1/2}$,
\begin{equation}\label{negregime}
\frac{\max_{i\in[n]}p_{ni}}{\spn}\build\to_{n\to\infty}^{}0. 
\end{equation}
Under this hypothesis, it has been shown \cite{almpimap02} that several features of the 
$\bp$-mapping, such as sizes of basins and number of cyclic points, can be described
asymptotically in terms of certain functionals of reflected Brownian bridge 
(this was originally proved in \cite{jpda94} for the uniform case 
$p_{ni}=1/n$). The two basic ingredients  % proving this are:
in the methodology of \cite{almpimap02} are:

\noindent (i)  Code the random mapping into a 
{\em mapping-walk} $H^{M_n}$ that contains enough 
information about the mapping; 

\noindent (ii) 
use a random bijection, called the Joyal correspondence \cite{joyal81}, 
that maps 
$\bp$-mappings into random doubly-rooted trees, called $\bp$-trees, whose
behavior is better understood. \\ In particular, 
the limits in law of associated encoding random walks can be shown to converge
to twice normalized Brownian excursion under condition (\ref{negregime}), and this information
lifts  back to  mappings, implying  that the  rescaled mapping  walks converge
weakly to twice standard reflecting Brownian bridge; that is, $\spn H^{M_n}\to 2\bb$
according to a certain topology on c\`adl\`ag functions. 
Results provable via this methodology encompass those proved in
\cite{op00} by somewhat different methods.

The goal of this paper is to extend this methodology to more general asymptotic
regimes for the distribution $\bp$, under the natural assumption
$\max_{i\in[n]}p_{ni}\to 0$  as $n\to  \infty$. % While encompassing  results of
%$\cite{op00}$ and the conjectured extensions therein, we mention that most our
%results are quite far from providing closed expressions. Other possible
%regimes besides (\ref{negregime}) are those for which
In these more general regimes,
several $\bp$-values are comparable to $\spn$ instead of being negligible.
Precisely, we will assume there exists $\btheta=(\theta_1,\theta_2,\ldots)$ such that 
\begin{equation}\label{regime}
\max_{i\in [n]}p_{ni}\build\to_{n\to\infty}^{}0 \qquad \mbox{ and }\qquad
\frac{p_{ni}}{\spn}\build\to_{n\to\infty}^{}\theta_i \, ,\qquad i\geq 1.
\end{equation}
By Fatou's Lemma, such a limiting $\btheta$ must satisfy $\sum_i\theta_i^2\leq
1$,   but  $\sum_i   \theta_i   $  may   be   finite  or   infinite.  We   let
$\theta_0=\sqrt{1-\sum_i\theta_i^2}$. 
A vertex $i\geq 1$ with $\theta_i>0$
then corresponds to a ``hub" \cite{jpda97ebac} or ``attracting center'' \cite{op00} for the mapping, because significantly
many more integers are likely to be  mapped to it as $n$ gets large than to those
for which $\theta_i=0$.  %(the notion of attracting center comes from
Our main result (Theorem \ref{T1}) roughly states that for $\bp_n$ satisfying
(\ref{regime}) with $\btheta=(\theta_1,\ldots,\theta_I,0,0,\ldots)$ with
$\theta_I>0$ and $\theta_0>0$ (the subset of such $\btheta$'s
is called 
$\Thfinite$), we have weak convergence 
\begin{equation}
\spn H^{M_n}\build\to_{}^{(d)} Z^{\sth}
\label{rough-limit}
\end{equation}
for a certain continuous process $Z^{\sth}$ to be described in section \ref{sec-TLP}, where the topology is in
general slightly weaker than the usual Skorokhod topology. We will also provide
criteria under which the stronger convergence  holds. In turn, we  will see how
this convergence and related results give information on the size of
the basins of  $M_n$, and on the  number of cyclic points, which  in the limit
arise as a kind of local time at $0$ for $Z^{\sth}$. 

To implement our methodology, the key point is that
under (\ref{regime}), the $\bp$-trees are known to converge in a certain sense 
(Proposition \ref{margptree}) to an 
{\em Inhomogeneous Continuum Random Tree} (ICRT) which we denote by
$\TT^{\sth}$. This family of trees was first investigated in \cite{jpda97ebac} % where
in the context of the additive coalescent.
What is important for this paper is the recent result
\cite{almpiep} that a certain class of ICRT's are encoded into
random excursion functions $H^{\sth}$, just as the Brownian tree is encoded into twice the
normalized Brownian excursion. % and we made these encoding functions $H^{\sth}$
The definition of $H^{\sth}$ is recalled in section \ref{sec-TLP},
where we also give the definition of 
the process
$Z^{\sth}$ as a functional of $H^{\sth}$. 

So the contribution of this paper is to show how the
ideas from
\cite{almpimap02} (in particular, the Joyal functional featuring in our Lemma \ref{lemj}) may be combined with the result of 
\cite{almpiep}
to prove the limit result indicated at (\ref{rough-limit}).
Once these ingredients are assembled, only a modest amount
of new technicalities
(e.g  part (ii) of Theorem \ref{thptree} and its use in the proof of
Theorem \ref{T1})
will be required.
One reason why ``only modest" is our restriction to the
case 
$\Thfinite$.
In \cite{almpiep} it is shown that the construction of $H^{\sth}$
and associated limit results for $\bp$-trees work in the
more general setting where $\sum_i \theta_i < \infty$.
It seems very likely that our new result (Theorem \ref{T1}) 
also extends to this setting, but the technicalities 
become more complicated.

While the existence of a limit process $Z^{\sth}$
provides qualitative information about aspects of the $\bp$-mappings,
enabling one to show that various limit distributions exist
and equal distributions of certain functionals of $Z^{\sth}$,
obtaining explicit formulas for such distributions remains
a challenging open problem.

\section{Statement of results}

\subsection{Mappings, trees, walks}
\label{sec-MTW}

We first introduce some % a certain number of
notation which is mostly taken from \cite{almpimap02}. 
If $m$ is a 
mapping on some finite set $S$, let ${\cal D}(m)$ be the directed graph with
vertex set $S$,
whose edges are $s\to m(s)$, and let $\cyc(m)$ be the set of cyclic
points, which is  further  partitioned  into disjoint  cycles,  $s$ and  $s'$
belonging to the same  cycle if one is mapped to the  other by some iterate of
$m$. For $c\in\cyc(m)$, if we remove the edges $c\to m(c)$ and $c'\to c$
where $c'$  is the unique point of  $S\cap{\cal C}(m)$ that is  mapped to $c$,
the component  of ${\cal D}(m)$ containing  $c$ is a tree  $\TT_c(m)$ which we
root at $c$. Label the disjoint cycles of $m$ as $\cyc_1(m),\cyc_2(m),\ldots$
with some ordering convention, then this in turn induces an
order on the {\em basins} of $m$: 
$$\com_j(m):=\bigcup_{c\in\cyc_j(m)}\TT_c(m).$$

\paragraph{$\bq$-biased order.}The ordering we will consider in this paper 
uses a convenient extra randomization, 
yet we mention that results similar to this
paper's could be established for different choices of basins ordering using
similar methods. See e.g.\
\cite{jpda02br}, where two different choices of ordering lead to two intricate
decompositions of Brownian bridge. Given $\bq$, a probability distribution on $S$ with
$q_s>0$ for every $s\in S$, consider an i.i.d.\ $\bq$-sample
$(X_2,X_3,\ldots)$ indexed by $\{2,3,\ldots\}$. If $m$ is a random mapping, we choose the $\bq$-sample
independently of $m$. Since $q_s>0$ for every $s\in S$, the following procedure
a.s.\ terminates:

$\bullet$ Let $\tau_1=2$ and let $\com_1(m)$ be the basin of $m$ containing
$X_2$. 

$\bullet$  If $\cup_{1\leq i\leq  j}\com_i(m)= S$  then end  the procedure;
else, given $\tau_j$ let $\tau_{j+1}=\inf\{k:X_k\notin \cup_{1\leq i\leq j}
\com_i(m)\}$ and let $\com_{j+1}$ be the basin containing $X_{\tau_{j+1}}$. 

This  induces  an order  on  basins  of $m$,  and  then  on the  corresponding
cycles. We add a further order on the cyclic points themselves by letting $c_j$
be the cyclic point of $\cyc_j(m)$ such that $X_{\tau_j}\in\TT_{c_j}$, and by
ordering the cyclic points within $\cyc_j(m)$ as follows:
$$m(c_j)\prec m^2(c_j)\prec\ldots\prec m^{|\cyc_j(m)|-1}(c_j)\prec c_j.$$
This extends to a linear order on $\cyc(m)$ by further letting $c_{j-1}\prec
m(c_j)$.  We call this (random) order on cyclic points and basins the
{\em $\bq$-biased random order}. In the special case where $\bq$ is the uniform
distribution on $S$, we call it the size-biased order. 

\paragraph{Coding trees and mappings with marked walks}

Let $\bT_n^o$ be the set of plane (ordered) rooted trees 
with $n$ labeled vertices $1,2,\ldots,n$, 
so that the children of any vertex $v$ are distinguished as
first, second, \ldots The cardinality  of $\bT_n^o$ is therefore $n!C_n$ where
$C_n$ is the $n$-th Catalan number. 
For any $T\in\bT_n^o$, we may put its set of vertices in
a special linear order $v_1,v_2,\ldots,v_n$ called {\em depth-first order}: 
we let  $v_1={\rm root}$, and  then $v_{j+1}$ is  the first (oldest)  child of
$v_j$ not in  $\{v_1,\ldots,v_j\}$ if any, or the oldest  brother of $v_j$ not
in $\{v_1,\ldots,v_j\}$ if  any, or the oldest brother of  the parent of $v_j$
not in $\{v_1,\ldots,v_j\}$, and so on. 
Write  $\HT^{T}(v)$ for  the  height of  vertex  $v$. For  any weight  sequence
$\bw=(w_1,\ldots,w_n)$ with $w_i>0$ for every $i$, let 
\begin{equation}\label{depth-first}
H^T_{\bw}(s)=\HT^T(v_i) \qquad \mbox{ if }\sum_{j=1}^{i-1}w_{v_j}\leq
s<\sum_{j=1}^iw_{v_j},
\end{equation}
and let $H^T_{\bw}(\sum_i w_i)=\HT^T(v_n)$. Call $H^T_{\bw}$ 
the height process of $T$. Notice that any $s\in[0,\sum_iw_i)$ 
specifies a vertex of $T$, which is $v_i$ in the case appearing in 
(\ref{depth-first}). We say that $v_i$ is {\em visited at time} $s$ by 
$H^T_{\bw}$. 
Intuitively, picture a particle touring the vertices in depth-first order
during the unit time interval, spending time $w_i$ at vertex $i$.

Given a mapping $m$ with basins and
cyclic points $c_1,\ldots,c_K$ in $\bq$-biased order for some $\bq$, we may
associate to each $\TT_{c_i}$ a walk as follows. First, turn these unordered
trees into plane trees by putting each set of children of each vertex 
in random exchangeable order, independently over vertices given $\TT_{c_i}$. 
Then associate to this ordered tree the height process $H^{\TT_{c_i}}_{\bw}$, 
with a
slight abuse of notation, where we are again given a weight function $\bw$ on
$[n]$ (though we use only the relevant labels appearing in $\TT_{c_i}$). We 
can now define the walk associated with $m$ to be
$$H^m_{\bw}(s)=H^{\TT_{c_i}}_{\bw}\left(s-\sum_{1\leq j<i}w(\TT_{c_j})\right) 
\qquad \mbox{ if
}\sum_{j=1}^{i-1}w(\TT_{c_j})\leq s<\sum_{j=1}^i w(\TT_{c_i}),$$
and $H^m_{\bw}(\sum_i w_i)=H^m_{\bw}(\sum_iw_i-)$, 
where $w(A)=\sum_{i\in A}w_i$. That
is, we  concatenate the tree-walks  associated to $\TT_{c_1},\ldots,\TT_{c_K}$
in this order. Again, there is a natural notion of vertex visited at time
$s<\sum_iw_i$. 

Further, let $D^m_{\bw}(i)=\sum_{j=1}^i w(\com_j(m))$ be the weight of the 
$i$ first basins, so that $D^m_{\bw}(i)$ is the time when the mapping-walk 
$H^m_{\bw}$ has completely
visited the vertices of the $i$-th basin, so $w(\com_i(m))=D^m_{\bw}(i)-
D^m_{\bw}(i-1)$ for
$i\geq 1$ with the convention  $D^m_{\bw}(0)=0$. We also let 
$\ell^m_{\bw}(s)$ be the number of
cyclic points that have been visited before time $s$, namely
$$\ell^m_{\bw}(s)=\sum_{j=1}^i{\bf 1}_{\left\{H^m_\bw (w(\{v_1,\ldots,v_j\})) =0\right\}}
\qquad\mbox{ whenever }\sum_{j=1}^{i-1}w_{v_j}\leq s<\sum_{j=1}^iw_{v_j},$$
with $\ell^m_{\bw}(\sum_iw_i)=\ell^m_{\bw}(\sum_iw_i-)$. 

\subsection{The Joyal functional}

We now define a functional $\bJ^u$ on the Skorokhod space $\D[0,1]$, which
translates into the world of encoding  paths the Joyal  bijection (recalled below) between trees
and mappings. Let $u\in[0,1]$.  Define the {\em pre-post
infimum} of $f\in\D[0,1]$ before and after $u$ to be the function 
$$s \to \un{f}_s(u)=\left\{\begin{array}{ll}
\inf_{t\in[s,u]}f_t & \mbox{ for }s<u\\
\inf_{t\in[u,s]}f_t & \mbox{ for }s\geq u.
\end{array}\right.$$
The function $\un{f}(u)$ is non-decreasing on $[0,u]$ and non-increasing on
$[u,1]$. If $[a,b]$ is a maximal flat interval for $\un{f}(u)$, we call the
recentered function $((f-\un{f}(u))(s+a),0\leq s\leq b-a)$ an excursion of $f$ 
above $\un{f}(u)$. Such a function may not be an excursion in the usual sense because it might
be zero for some $s\in(0,b-a)$.  Further, if two distinct such intervals
$[a,b]$ and $[c,d]$ satisfy $f(b)=f(c)$, then it must be that $b<u<c$, and in
this case we call the function obtained by concatenating the excursion of $f$
above $\un{f}(u)$ on $[a,b]$ and $[c,d]$ a (generalized) excursion of $f$
above $\un{f}(u)$. Label as $\eps_1,\eps_2,\ldots$ the generalized 
excursions of $f$ above $\un{f}(u)$, according to decreasing durations
$l_1,l_2,\ldots$.  Write also $h_i$ for the ``height'' of the excursion
$\eps_i$, i.e.\ the value taken by $\un{f}(u)$ on the flat interval of the
excursion. We  define a function  $\bJ^u(f)$ that arranges these  excursions in
order of heights:
\begin{equation}\label{joyalfunct}
\bJ^u(f)(s)=\eps_i\left(s-\sum_{j:h_j<h_i}l_j\right)\qquad 
\mbox{ if }\sum_{j:h_j<h_i}l_j\leq s<\sum_{j:h_j\leq h_i}l_j,
\end{equation}
with the convention that $\bJ^u(f)(s)=0$ on $[\sum_il_i,1]$. 

To keep track of the structure of the original function, we finally add marks
at the points $g^u_i(f)=\sum_{j:h_j<h_i}l_j$ and $d^u_i(f)=\sum_{j:h_j\leq
h_i}l_j,i\geq 1$.  In particular, if $\bJ^u(f)$ if non-zero on
$(g^u_i(f),d^u_i(f))$, then the  $\eps_i$ is an ``usual'' excursion 
rather than ``generalized'' excursion. 

\subsection{The limiting process and main result}
\label{sec-TLP}

Let us recall the construction \cite{almpiep} of 
the exploration process of the ICRT $\TTth$
for $\btheta\in\Thfinite$.  Let $(b_s,0\leq s\leq 1)$ be a standard Brownian
bridge, $U_1,\ldots,U_I$ be independent uniform random variables independent
of $b$, and 
$$\xbrt(s)=\theta_0 b_s+\sum_{i=1}^I\theta_i({\bf 1}_{\{U_i\leq s\}}-s)\, , \qquad
0\leq s\leq 1.$$
Such a process  has a.s.\ a unique time where it  attains its overall minimum,
and this  time is a  continuity time, call  it $s_{\min}$. Define  the {\em Vervaat
transform} of $\xbrt$ by 
$$X^{\sth}(s)=\xbrt(s+s_{\min} [\mbox{mod }\, 1])-\xbrt(s_{\min})\, ,\qquad 
0\leq s\leq 1,$$
and let $t_i=U_i-s_{\min}[\mbox{mod }\, 1],1\leq i\leq I$ be the jump times of
$X^{\sth}$. Let $T_i=\inf\{s\geq t_i:X^{\sth}_s=X^{\sth}_{t_i-}\}$ and write 
$$R^{\sth}_i(s)=\left\{\begin{array}{cl}
\inf_{t_i\leq u\leq s}X^{\sth}_u -X^{\sth}_{t_i-}& \mbox{ if }s\in[t_i,T_i]\\
0& \mbox{ else. }
\end{array}\right. 
$$
Last,   let   $Y^{\sth}=X^{\sth}-\sum_{i=1}^IR^{\sth}_i$.   This  process 
$Y^{\sth}$ is
continuous, and it  is intuitively described by: ``take away  all the jumps of
$X^{\sth}$ and reflect  the process above its infimum  after these jumps until
$X^{\sth}$ gets back to the level it started at before jumping''. 
The exploration process of $\TT^{\sth}$ is then defined as
$H^{\sth}=\sfrac{2}{\theta_0^2}Y^{\sth}$. 

The open set  $\{s\in(t_i,T_i):H^{\sth}(s)>H^{\sth}(t_i)\}$ associated with jump $i$ can be decomposed into
disjoint open  intervals $(t_{ij},T_{ij}),j\geq 1$, ranked  by decreasing order
of lengths. 

Now take a  uniform$(0,1)$ random variable $U$ independent  of $Y^{\sth}$, and
consider the process $Z^{\sth}=\bJ^U(H^{\sth})$.  Recall that this process has
marks $g_i^U(H^{\sth}),d_i^U(H^{\sth})$, which we more simply call $g_i,d_i$. 
The following facts are
simple consequences of usual properties of Brownian motion and Brownian
bridge: 
\begin{itemize}
\item The sum of durations of generalized excursions of $H^{\sth}$ above
$\un{H}^{\sth}(U)$ is $1$, meaning $\sum_{i\geq 1}(d_i-g_i)=1$. 
\item The corresponding excursion 
heights $h_i,i\geq 1$ are a.s.\ everywhere dense in
$[0,H^{\sth}_U]$. 
\end{itemize}
Now let $V_1,V_2,\ldots$  be independent uniform$(0,1)$ variables, independent
of $U$ and $H^{\sth}$.  Define recursively a
sequence $D_0=0<D_1<D_2<\ldots< 1$ by
$$D_n=\inf\{s: s>D_{n-1}+V_n(1-D_{n-1})\mbox{ and }
\exists i\geq 1,s=d_i\}\qquad n\geq 1.$$ 
Last, we define the local time function of $Z^{\sth}$ as follows: for $s$
in an excursion interval of $Z^{\sth}$ above $0$, let $L^{\sth}_s$ be the
``height'' of the corresponding generalized excursion of $H^{\sth}$ above
$\un{H}^{\sth}$.  This defines $L^{\sth}$ on a dense subset of $[0,1]$ 
as an increasing function, which can be extended to the whole 
interval $[0,1]$ uniquely as a continuous function, because $H^{\sth}$ is
itself continuous, and the excursion heights are dense in
$[0,H^{\sth}_U]$. Notice that this ``local'' time has the unusual property
that its increase times do not exactly match with the zero set of
$Z^{\sth}$; rather, the set of increase times is the closure of
$\{g_i,d_i,i\geq 1\}$. 

Now let $\bp_n,\bq_n,\bw_n$ be three sequences of probabilities on 
$[n]$ charging every point. 
Consider a $\bp_n$-mapping $M_n$ with basins in $\bq_n$-biased order, and let
$H^{M_n}_{\bw_n}$ be the associated walk. We let $H^{M_n}:=H^{M_n}_{\bp_n}$.
Our main result is 
\begin{thm}\label{T1}
Suppose
 $\max_i q_{ni}\to 0$ as $n\to \infty$. 

{\rm (i)} Under the asymptotic regime (\ref{regime}) for $\bp_n$, with 
limiting $\btheta\in\Thfinite$, and if (\ref{extrahyp1},\ref{extrahyp2}) 
below are satisfied, then for any weight function $\bw$ satisfying 
$\max_iw_{ni}\to 0$, we have the convergence in law in the usual Skorokhod
topology on $\D[0,1]$
$$\spn H^{M_n}_{\bw_n}\build\to_{}^{(d)} Z^{\sth}.$$

{\rm (ii)} Moreover, jointly with the above convergence, the marks 
$D^{M_n}_{\bw_n}(1),D^{M_n}_{\bw}(2),\ldots$ converge in law to 
$D_1,D_2,\ldots$. 

{\rm (iii)} Jointly with the above convergences, 
$\spn\ell^{M_n}_{\bw_n}\build\to_{}^{(d)} L^{\sth}$ for the uniform topology. 

{\rm (iv)} In general, 
under the asymptotic regime (\ref{regime}) for $\bp_n$, with 
limiting $\btheta\in\Thfinite$, 
one has convergence in law for the $*$-topology defined  
in \cite{jpda02inv}
$$\spn H^{M_n}\build\to_{}^{(d)} Z^{\sth},$$
and the convergences of {\rm (ii),(iii)} hold jointly for $\bw_n=\bp_n$. 
\end{thm}

We echo \cite[Corollary 1]{almpimap02} by stating

\begin{crl}\label{C1}
Under (\ref{regime},\ref{extrahyp1},\ref{extrahyp2}) 
with finite-length limiting
$\btheta$, and for any weight function $\bw_n,\bq_n$ with
$\max_i \max(w_{ni},q_{ni})\to0$ as $n\to\infty$, we have 
$$(\bw_n(\com_j(M_n)),\spn\card\cyc_j(M_n),j\geq 1)\build\to_{}^{(d)}
(D_j-D_{j-1},L^{\sth}_{D_j}-L^{\sth}_{D_{j-1}},j\geq 1).$$
\end{crl}
Notice that for uniform $\bw_n$, the first component equals 
$n^{-1}\card(\com_j(M_n))$.  

The essential point of the $*$-topology is the following
property \cite{almpimap02}.
One has $f^n \to f$ 
for the $*$-topology, where $f_n \in\D[0,1]$ and $f\in C[0,1]$ is continuous,
if and only of there exist $g_n,h_n\in\D[0,1]$ with $f_n=g_n+h_n$
such that $g_n\to f$ uniformly, $h_n\geq 0$ and ${\rm
Leb}\{x:h_n(x)>0\}\to0$. % In particular, Theorem \ref{T1} (iv) allows to deduce
Thus the $*$-convergence asserted in (iv) is compatible with
the possible presence
of upward
``spikes'' on the mapping-walk, which have 
arbitrary large height but vanishing weight. % One has $f^n\to f$
In particular, Theorem \ref{T1} (iv) allows us to deduce
the asymptotic ``height'' (distance to the set $\cyc(M_n)$) of a randomly
$\bp_n$-chosen vertex, but not the behavior of the asymptotic
maximum height over all vertices, which is however handled under the hypotheses in
(i). Under the same hypotheses, we can handle quantities such as the diameter
of the random mapping (the maximal $k$ such that there exists $v$ with
$v,m(v),\ldots,m^{k-1}(v)$ pairwise distinct). 

Although this result leaves a large degree of freedom for choosing the order of
basins, we stress that other  orderings are possible, such as ordering the basins
according to increasing order of the least vertices they contain, or ordering
cycles by order of least vertex they contain. The first order is in fact
equivalent to the size-biased order described above, up to relabeling, and
the second order could be also handled by our methods, although the marks
$D_i$ would have to be defined in a different way, see \cite{jpda02br}. 

Last, we stress that the hypotheses (\ref{extrahyp1},\ref{extrahyp2})
below are by no means necessary, we believe that they are in fact quite
crude (see \cite{almpiep} for further discussion). Also, as  discussed below, we believe that Theorem
\ref{T1} (iv) remains true for much more general $\bw_n$. 

\section{Proofs}

\subsection{p-trees and associated walks}

\paragraph{$\bp$-trees and their walks.}
We now define the random trees whose asymptotics are related to the
process $H^{\sth}$, namely $\bp$-trees. Let $\bT_n$ be the set constituted of
the $n^{n-1}$ (unordered) rooted labeled trees on $[n]$.  For $\bp$ a
probability measure charging every point of $[n]$, let $\TT^{\bp}$ be the
random variable in $\bT_n$ with law 
\begin{equation}\label{ptreedef}
P(\TT^{\bp}=t)=\prod_{i\in[n]}p_i^{c_i(t)}\, ,\qquad t\in\bT_n,
\end{equation}
where $c_i(t)$ is the number of children of $i$ in $t$. The fact that 
(\ref{ptreedef}) indeed defines a {\em probability} measure amounts to the 
Cayley  multinomial expansion for  trees \cite{jp01hur}.  For $t\in\bT_n$,  we can
associate a random $\bT_n^o$-valued tree $t^o$ 
by putting each set of children of a
given vertex in uniform random order, independently over distinct vertices, so
given a weight function $\bw$ on $[n]$ we may associate to $t$ the random walk 
$H^t:=H^{t^o}$ as defined in section \ref{sec-MTW}.  We will now apply this to the random trees
$\TT^{\bp}$ and their associated height processes $H^{\bp}_{\bw}:=
H^{\TT^{\bp}}_{\bw}$. When $\bw=\bp$ we let $H^{\bp}:=H^{\bp}_{\bp}$. 

\paragraph{Asymptotics.}
We introduce two extra hypothesis on the sequence $\bp_n$ besides
(\ref{regime}). The first one prevents exponentially small (in the scale $\sp$)
$\bp$-values from appearing:
\begin{equation}\label{extrahyp1}
(\min_i p_{ni})^{-1}=o(\exp(\alpha/\spn))\, ,\qquad \forall\, \alpha>0.
\end{equation}
The second states that ``small'' $\bp$-values are of rough order
$\sp^2$. Suppose there exists some non-negative finite r.v.\ $Q$
such that, letting $\bar{\bp}=(0,\ldots,0,p_{I+1},p_{I+2},\ldots)$, 
\begin{equation}\label{extrahyp2}
\lim_{n\to\infty}E\left[\exp\left(\frac{\lambda \bar{p}_{n\xi}}{\spn^2}
\right)\right]=E[\exp(\lambda Q)]<\infty
\end{equation}
for every $\lambda$ in a neighborhood of $0$.  Here, $\xi$ denotes a random
variable with law $\bp$, so $\bar{p}_{n\xi}$ is its $\bar{\bp}_n$-value. 

The key results on $\bp$-trees are the following
variations of \cite[Theorems 1,3]{almpiep}. For $k\geq 1$, let
$X_2,\ldots,X_k$ be independent $\bp$-sampled vertices of $\TT^{\bp}$,
independent of $\TT^{\bp}$. Let $r_k(\TT^{\bp})$ be the subtree of $\TT^{\bp}$
spanned by the root and $X_2,\ldots,X_k$, re-interpreted as a {\em tree with
edge-lengths}, in the sense that two vertices separated by a single edge are at
distance $1$, and we delete all the nodes that have degree $2$, so the distance
between two vertices on the final tree is equal to the number of deleted nodes
plus $1$. The tree
$r_k(\TT^{\bp})$ is thus a discrete rooted tree with at most $k$ leaves, which
has no degree $2$ vertices, and with
lengths attached to each of its edges. The notion of convergence on the space
of trees with edge-length is the usual convergence for the product topology, so
$\bt_n\to \bt$ if both trees have the same discrete structures for all sufficiently large
$n$, and the vector of edge-lengths of $\bt_n$ converges to that of $\bt$. Last,
for $a>0$ we let $a\otimes \bt$ be the tree with edge-length with same discrete
structure as $\bt$, and where all distances have been multiplied by $a$. 
\begin{prp}[\cite{jpmc97b}]\label{margptree}
Suppose that $\bp_n$ satisfies (\ref{regime}).
For every $k$, the tree $\spn\otimes r_k(\TT^{\bp_n})$ converges in distribution
to $\TT^{\sth}_k$, the $k$-th marginal of  the ICRT described in
section \ref{sec-ICRT}. 
\end{prp}

\begin{thm}\label{thptree}
{\rm (i)}
Suppose that $\bp_n$ satisfies (\ref{regime},\ref{extrahyp1},\ref{extrahyp2}), 
with $\btheta\in\Thfinite$, and that $\bw_n$ satisfies
$\max_iw_{ni}\to0$. Then 
$$\spn H^{\bp_n}_{\bw_n}\build\to_{}^{(d)} H^{\sth}$$
for the usual Skorokhod topology (and hence for the uniform topology since the limit
is continuous).

{\rm (ii)}
Under the assumptions of (i), 
for each $1\leq i\leq I$, there exist random sequences $t^{\bp_n}_i,T^{\bp_n}_i$
and $t^{\bp_n}_{ij},T^{\bp_n}_{ij},j\geq 1$ with
$$H^{\bp_n}_{\bw_n}(t^{\bp_n}_i)=
H^{\bp_n}_{\bw_n}(t^{\bp_n}_{ij})=H^{\bp_n}_{\bw_n}
(T^{\bp_n}_{ij})\mbox{ for every } j\geq 1,$$ 
and $H^{\bp_n}_{\bw_n}(s)\geq H^{\bp_n}_{\bw_n}(t^{\bp_n}_i)$ for
$s\in[t^{\bp_n}_i,T^{\bp_n}_i]$, such that jointly with the convergence of {\rm
(i)}, one has convergence in law
$$(t^{\bp_n}_i,T^{\bp_n}_i,t^{\bp_n}_{ij},T^{\bp_n}_{ij},1\leq i\leq I,j\geq 1)
\build\to_{}^{(d)} (t_i,T_i,t_{ij},T_{ij},1\leq i\leq I,j\geq 1),$$
with the notations of section \ref{sec-TLP}. 

{\rm (iii)} 
Suppose only that $\bp_n$ satisfies (\ref{regime}), then
$$\spn H^{\bp_n}\build\to_{}^{(d)} H^{\sth}$$
in the $*$-topology. Moreover, the statement of {\rm (ii)} still holds
for $\bw_n=\bp_n$. 
\end{thm}

\proof
Except for the last sentence, 
point (iii) is a consequence of \cite[Proposition 7]{jpda02inv} which states that the
convergence of marginals of $\bp$-trees to that of the limiting ICRT
\cite[Proposition 1 and (23)]{almpiep} is equivalent to the $*$-convergence of
the rescaled walk $\spn H^{\bp_n}$ with weights $\bw_n=\bp_n$ to 
$H^{\sth}$. 

Point (i) was proved in \cite[Theorem 3, Corollary 3]{almpiep} in the two special
cases  where $\bw_n=\bp_n$ and where $\bw_n=(1/n,\ldots,1/n)$  ($n$  times).  The
general case uses the same proof  as Corollary 3 in the stated paper.
%, starting
%from the convergence with weights $(1/n,\ldots,1/n)$ (call $H^n$ the
%associated tree-walk). 
By the weak law of
large numbers for sampling without replacement applied to $\bw_n$, we have 
$\sup_{0\leq t\leq 1}|S_{\bw_n,0}(t)-t|\to 0$ in probability, where 
$S_{\bw_n,0}$ is the linear interpolation between the points $((\sum_{1\leq k\leq
i}w_{n\pi(k)},i/n),1\leq i\leq n)$, and  where $\pi$ is a uniformly distributed
random permutation on $[n]$.  This implies the result because, as shown in
\cite{almpiep}, the depth-first order on vertices $v_1,v_2,\ldots,v_n$
of  a $\bp_n$-tree  is a  (random) shift  of a  uniform permutation  of $[n]$.
Therefore, the linear interpolation $S_{\bw_n}$ between points 
$((\sum_{1\leq k\leq i}w_{nv_k},i/n),1\leq i\leq n)$ also uniformly converges
to the identity, and the conclusion follows from the fact that $H^{\bp_n}_{\bw_n}=H^{\bp_n}
\circ(S_{\bp_n})^{-1}\circ S_{\bw_n}$. 

\setlength{\unitlength}{0.025in}
\begin{picture}(250,125)(0,-50)
\multiput(0,40)(1,0){100}{\line(1,0){0.5}}
\put(10,40){\line(0,1){2}}
\put(10,42){\line(1,0){1}}
\put(11,42){\line(0,1){3}}
\put(11,45){\line(1,0){1}}
\put(12,45){\line(0,-1){1}}
\put(12,44){\line(1,0){1}}
\put(13,44){\line(0,1){3}}
\put(13,47){\line(1,0){1}}
\put(14,47){\line(0,-1){1}}
\put(14,46){\line(1,0){1}}
\put(15,46){\line(0,1){3}}
\put(15,49){\line(1,0){1}}
\put(16,49){\line(0,-1){1}}
\put(16,48){\line(1,0){1}}
\put(17,48){\line(0,-1){2}}
\put(17,46){\line(1,0){1}}
\put(18,46){\line(0,-1){2}}
\put(18,44){\line(1,0){1}}
\put(19,44){\line(0,1){1}}
\put(19,45){\line(1,0){1}}
\put(20,45){\line(0,-1){3}}
\put(20,42){\line(1,0){1}}
\put(21,42){\line(0,1){1}}
\put(21,43){\line(1,0){1}}
\put(22,43){\line(0,-1){2}}
\put(22,41){\line(1,0){1}}
\put(23,41){\line(0,-1){1}}
\put(40,40){\line(0,1){2}}
\put(40,42){\line(1,0){1}}
\put(41,42){\line(0,1){3}}
\put(41,45){\line(1,0){1}}
\put(42,45){\line(0,-1){1}}
\put(42,44){\line(1,0){1}}
\put(43,44){\line(0,1){3}}
\put(43,47){\line(1,0){1}}
\put(44,47){\line(0,-1){2}}
\put(44,45){\line(1,0){1}}
\put(45,45){\line(0,1){2}}
\put(45,47){\line(1,0){1}}
\put(46,47){\line(0,1){3}}
\put(46,50){\line(1,0){1}}
\put(47,50){\line(0,-1){1}}
\put(47,49){\line(1,0){1}}
\put(48,49){\line(0,1){3}}
\put(48,52){\line(1,0){1}}
\put(49,52){\line(0,-1){2}}
\put(49,50){\line(1,0){1}}
\put(50,50){\line(0,-1){3}}
\put(50,47){\line(1,0){1}}
\put(51,47){\line(0,1){1}}
\put(51,48){\line(1,0){1}}
\put(52,48){\line(0,-1){2}}
\put(52,46){\line(1,0){1}}
\put(53,46){\line(0,-1){1}}
\put(53,45){\line(1,0){1}}
\put(54,45){\line(0,-1){3}}
\put(54,42){\line(1,0){1}}
\put(55,42){\line(0,1){1}}
\put(55,43){\line(1,0){1}}
\put(56,43){\line(0,-1){2}}
\put(56,41){\line(1,0){1}}
\put(57,41){\line(0,-1){1}}
\put(75,40){\line(0,1){2}}
\put(75,42){\line(1,0){1}}
\put(76,42){\line(0,1){3}}
\put(76,45){\line(1,0){1}}
\put(77,45){\line(0,-1){1}}
\put(77,44){\line(1,0){1}}
\put(78,44){\line(0,1){3}}
\put(78,47){\line(1,0){1}}
\put(79,47){\line(0,1){2}}
\put(79,49){\line(1,0){1}}
\put(80,49){\line(0,-1){2}}
\put(80,47){\line(1,0){1}}
\put(81,47){\line(0,1){1}}
\put(81,48){\line(1,0){1}}
\put(82,48){\line(0,-1){2}}
\put(82,46){\line(1,0){1}}
\put(83,46){\line(0,-1){1}}
\put(83,45){\line(1,0){1}}
\put(84,45){\line(0,-1){3}}
\put(84,42){\line(1,0){1}}
\put(85,42){\line(0,1){1}}
\put(85,43){\line(1,0){1}}
\put(86,43){\line(0,-1){2}}
\put(86,41){\line(1,0){1}}
\put(87,41){\line(0,-1){1}}
\put(10,0){\line(0,1){2}}
\put(10,2){\line(1,0){1}}
\put(11,2){\line(0,1){3}}
\put(11,5){\line(1,0){1}}
\put(12,5){\line(0,-1){1}}
\put(12,4){\line(1,0){1}}
\put(13,4){\line(0,1){3}}
\put(13,7){\line(1,0){1}}
\put(14,7){\line(0,-1){1}}
\put(14,6){\line(1,0){1}}
\put(15,6){\line(0,1){3}}
\put(15,9){\line(1,0){1}}
\put(16,9){\line(0,-1){1}}
\put(16,8){\line(1,0){1}}
\put(17,8){\line(0,-1){2}}
\put(17,6){\line(1,0){1}}
\put(18,6){\line(0,-1){2}}
\put(18,4){\line(1,0){1}}
\put(19,4){\line(0,1){1}}
\put(19,5){\line(1,0){1}}
\put(20,5){\line(0,-1){3}}
\put(20,2){\line(1,0){1}}
\put(21,2){\line(0,1){1}}
\put(21,3){\line(1,0){1}}
\put(22,3){\line(0,-1){2}}
\put(22,1){\line(1,0){1}}
\put(23,1){\line(0,-1){1}}
\put(31,-2){\line(1,0){1}}
\put(32,-2){\line(0,1){2}}
\put(32,0){\line(1,0){1}}
\put(31,-2){\line(0,1){3}}
\put(31,1){\line(-1,0){1}}
\put(40,0){\line(0,1){2}}
\put(40,2){\line(1,0){1}}
\put(41,2){\line(0,1){3}}
\put(41,5){\line(1,0){1}}
\put(42,5){\line(0,-1){1}}
\put(42,4){\line(1,0){1}}
\put(43,4){\line(0,1){3}}
\put(43,7){\line(1,0){1}}
\put(44,7){\line(0,-1){2}}
\put(44,5){\line(1,0){1}}
\put(45,5){\line(0,1){2}}
\put(45,7){\line(1,0){1}}
\put(46,7){\line(0,1){3}}
\put(46,10){\line(1,0){1}}
\put(47,10){\line(0,-1){1}}
\put(47,9){\line(1,0){1}}
\put(48,9){\line(0,1){3}}
\put(48,12){\line(1,0){1}}
\put(49,12){\line(0,-1){2}}
\put(49,10){\line(1,0){1}}
\put(50,10){\line(0,-1){3}}
\put(50,7){\line(1,0){1}}
\put(51,7){\line(0,1){1}}
\put(51,8){\line(1,0){1}}
\put(52,8){\line(0,-1){2}}
\put(52,6){\line(1,0){1}}
\put(53,6){\line(0,-1){1}}
\put(53,5){\line(1,0){1}}
\put(54,5){\line(0,-1){3}}
\put(54,2){\line(1,0){1}}
\put(55,2){\line(0,1){1}}
\put(55,3){\line(1,0){1}}
\put(56,3){\line(0,-1){2}}
\put(56,1){\line(1,0){1}}
\put(57,1){\line(0,-1){1}}
\put(66,-7){\line(1,0){1}}
\put(67,-7){\line(0,1){2}}
\put(67,-5){\line(1,0){1}}
\put(66,-7){\line(0,1){3}}
\put(66,-4){\line(-1,0){1}}
\put(75,0){\line(0,1){2}}
\put(75,2){\line(1,0){1}}
\put(76,2){\line(0,1){3}}
\put(76,5){\line(1,0){1}}
\put(77,5){\line(0,-1){1}}
\put(77,4){\line(1,0){1}}
\put(78,4){\line(0,1){3}}
\put(78,7){\line(1,0){1}}
\put(79,7){\line(0,1){2}}
\put(79,9){\line(1,0){1}}
\put(80,9){\line(0,-1){2}}
\put(80,7){\line(1,0){1}}
\put(81,7){\line(0,1){1}}
\put(81,8){\line(1,0){1}}
\put(82,8){\line(0,-1){2}}
\put(82,6){\line(1,0){1}}
\put(83,6){\line(0,-1){1}}
\put(83,5){\line(1,0){1}}
\put(84,5){\line(0,-1){3}}
\put(84,2){\line(1,0){1}}
\put(85,2){\line(0,1){1}}
\put(85,3){\line(1,0){1}}
\put(86,3){\line(0,-1){2}}
\put(86,1){\line(1,0){1}}
\put(87,1){\line(0,-1){1}}

\put(-15,30){\line(0,1){25}}
\put(-15,40){\line(1,0){2}}
\put(-21,38){$h$}
\put(-15,-20){\line(0,1){35}}
\put(-15,0){\line(1,0){2}}
\put(-21,-2){$h$}
\put(-15,-2){\line(1,0){2}}
\put(-11,-4){$\scriptstyle{h_n(1)}$}
\put(-15,-7){\line(1,0){2}}
\put(-28,-9){$\scriptstyle{h_n(2)}$}
\put(0,-40){\line(1,0){100}}
\put(0,-42){\line(0,1){4}}
\put(-2,-48){$t_i$}
\put(10,-42){\line(0,1){4}}
\put(8,-48){$t_{ij(1)}$}
\put(23,-42){\line(0,1){4}}
\put(21,-48){$T_{ij(1)}$}
\put(40,-42){\line(0,1){4}}
\put(38,-48){$t_{ij(2)}$}
\put(57,-42){\line(0,1){4}}
\put(55,-48){$T_{ij(2)}$}
\put(75,-42){\line(0,1){4}}
\put(73,-48){$t_{ij(3)}$}
\put(87,-42){\line(0,1){4}}
\put(85,-48){$T_{ij(3)}$}
\put(100,-42){\line(0,1){4}}
\put(98,-48){$T_i$}
\put(13.5,-42){\line(0,1){4}}
\put(11.5,-35){$U_1$}
\put(48.5,-42){\line(0,1){4}}
\put(46.5,-35){$U_2$}
\put(81.5,-42){\line(0,1){4}}
\put(79.5,-35){$U_3$}
\put(135,20){\line(0,1){35}}
\put(135,40){\line(1,0){2}}
\put(129,38){$h$}
\put(135,38){\line(1,0){2}}
\put(139,36){$\scriptstyle{h_n(1)}$}
\put(135,33){\line(1,0){2}}
\put(122,31){$\scriptstyle{h_n(2)}$}
\put(175,52){\line(0,-1){62}}
\put(175,40){\line(-1,0){7}}
\put(175,40){\line(1,0){8}}
\put(215,52){\line(0,-1){62}}
\put(215,38){\line(-1,0){9}}
\put(215,33){\line(1,0){15}}
\put(173,54){$2$}
\put(213,54){$2$}
\put(163,38){$1$}
\put(202,36){$1$}
\put(185,38){$3$}
\put(232,31){$3$}
\put(169,-16){root}
\put(209,-16){root}
\put(180,-45){\bf Figure 1.}
\end{picture}

%\vspace{0.3in}
%{\bf Figure 1.}

%\vspace{0.2in}
\noindent
Point (ii) refines one aspect of (i).  First consider the case $\bw_n = \bp_n$.   By Skorokhod's representation
theorem, suppose that the convergence in law of (i) holds almost-surely.
Fix $i$.
Figure 1 shows schematically (top left) three of the excursions
of $H^{\sth}$ associated with jump $i$.
All have the same height, $h$ say.
The lower left diagram in Figure 1 shows corresponding parts
of $H^{\bp_n}$.
Consider
the minimum value $h_n(1)$ of $H^{\bp_n}$ between $T_{ij(1)}$ and $t_{ij(2)}$,
and
the minimum value $h_n(2)$ of $H^{\bp_n}$ between $T_{ij(2)}$ and $t_{ij(3)}$.
The key claim is
\begin{equation}
h_n(1) = h_n(2) \mbox{ for all large } n .
\label{claim}
\end{equation}
To verify (\ref{claim}),
%Fix $i$. By uniform convergence, we may find $t^{\bp_n}_i,T^{\bp_n}_i$ 
%converging to $t_i,T_i$ almost-surely, satisfying $H^{\bp_n}\geq
%H^{\bp_n}(t^{\bp_n}_i)$ on
%$[t^{\bp_n}_i,T^{\bp_n}_i]$.  Since $H^{\bp_n}$ increases only by jumps with
%magnitude $1$, with no loss of generality, we may also
%choose these sequences such that there exists $m^{\bp_n}_1$ with
%$H^{\bp_n}(m^{\bp_n}_1)=H^{\bp_n}(t^{\bp_n}_i)$ converging to 
%some
%$m_1\in(t_i,T_i)$ with $H^{\sth}(m_1)=H^{\sth}(t_i)$, so $m_1$ is in the 
%closure of the set $\{t_{ij},T_{ij},j\geq 1\}$. Then if (ii) did not hold, we 
%would find some other local minimum $m^{\bp_n}_2$ of $H^{\bp_n}$ 
%converging to some $m_2\in(t_i,T_i)$ with $H^{\sth}(m_2)=H^{\sth}(m_1)$, but 
%with $H^{\bp_n}(s)\geq 
%H^{\bp_n}(m^{\bp_n}_2)>H^{\bp_n}(m^{\bp_n}_1)$ for every $s$ in a (fixed) 
%neighborhood of $m_2$ and every $n$ large. 
%Suppose $m_1<m_2$ and take three independent uniform random variables 
take three independent uniform random variables
$U_1,U_2,U_3$ on $[0,1]$ independent of $H^{\sth},H^{\bp_n},n\geq 1$. These 
random  variables  specify  three  $\bp_n$-chosen vertices  on  $\TT^{\bp_n}$,
namely those which are visited by $H^{\bp_n}$ at these times. %Then 
On an event of positive probability
we have $U_k \in (t_{ij(k)},T_{ij(k)}), \ k = 1,2,3$.
%the probability that: $U_1<U_2<U_3, 
%\inf_{u\in[U_1,U_2]}H^{\bp_n}(u)=H^{\bp_n}(m^{\bp_n}_1)$ 
%and $\inf_{u\in[U_2,U_3]}H^{\bp_n}(u)=H^{\bp_n}(m^{\bp_n}_2)$ for
%every $n$ large enough 
%is positive. This proves that the subtree of $\TT^{\bp_n}$
Consider the subtree of $\TT^{\bp_n}$
spanned by the root and the three vertices encoded by $U_1,U_2,U_3$.  %has one of
If $h_n(1) \neq h_n(2)$ then, on the above event, the
subtree has an edge of length 
$|h_n(2)-h_n(1)|$ 
(as shown in rightmost tree in Figure 1),
but this is not converging to the correct limit asserted
in Proposition \ref{margptree} 
(in the sense of convergence of discrete structures mentioned
above Proposition \ref{margptree})
because the limit tree (the second-right tree in Figure 1)
has different tree shape.
Thus we can deduce (\ref{claim}) using Proposition \ref{margptree}.
It is then straightforward to deduce the full assertion of (ii)
from the case (\ref{claim}) of three excursions.

%its edges with length
%$m^{\bp_n}_2-m^{\bp_n}_1>0$. 
%When scaling this tree by $\spn$, the length of the
%corresponding vertex tends to $0$ with positive probability, but it is never
%$0$, so the discrete structure of $r_3(\TT^{\bp_n})$ can by no mean be
%stationary on this event, contradicting Proposition \ref{margptree}. 

Treating the case of general weights $\bw_n$ is done by asking (again by the
Skorokhod representation theorem) that the uniform convergence of 
$S^{-1}_{\bw_n}\circ S_{\bp_n}$ to identity is also almost-sure. Then replace 
$U_k,k=1,2,3$ by  $U_k^{\bw_n}=S^{-1}_{\bw_n}\circ S_{\bp_n}(U_k),k=1,2,3$, so
the new variables encode again $\bp_n$-chosen vertices.  %and
%ask that $U^{\bw_n}_1<U^{\bw_n}_2<U^{\bw_n}_3$ only for $n$ large enough, which
%holds with positive probability. The case of $*$-convergence (for 
The case of $*$-convergence (for 
$\bw_n=\bp_n$) is similar (see also the proof of \cite[Lemma 2]{almpimap02}). 
\cq

\rem
To prove (iii) for more general weights $\bw_n$, we
could try to use the same method as above (first treating the case of uniform
weights). But if $f_n\to f$ for the
$*$-topology with $f$ continuous, and if $S_n$ is a strictly increasing
piecewise linear continuous function that converges uniformly to the
identity on $[0,1]$, then $f_n\circ S_n$ need not converge to $f$ for the
$*$-topology. Indeed, with the above notation, this convergence is equivalent
to ${\rm Leb}\{x\in[0,1]:h_n\circ S_n(x)>0\}\to0$. But this last quantity is 
$\int_0^1 {\bf 1}_{\{h_n>0\}}(S^{-1}_n)'(x)\d x$. So we would need a sharper
result than the weak law of large numbers for sampling without replacemement
to estimate the values of the derivative at points where $h_n>0$.  
% For example, we get the result in the
%somehow trivial case where the $w_{ni}$'s are all of order $\Theta(1/n)$ and 
%for $i>I$, the $p_{ni}$'s are all of order
%$\Theta(1/n)$, so the derivatives of associated interpolations and their 
%inverses are bounded away from $0$ and
%$\infty$ (except maybe for $\bp_n$ on a finite number of intervals with 
%uniformly
%vanishing size of order $1/\sqrt{n}$). However, it was proved in
However, it was proved in 
\cite[Theorem 25]{jpda02inv} by different methods 
that in the asymptotically negligible regime 
(\ref{negregime}), Theorem \ref{T1} (i) is still valid for general weights
$\bw_n$ satisfying $\max_iw_{ni}\to0$. It would therefore be surprising
if the same result did not hold here. 

\subsection{The Joyal correspondence.}

Let us now describe the Joyal correspondence between trees and mappings,
designed to push the distribution of $\bp$-trees onto the distribution of
$\bp$-mappings. Let $\bq$ be a probability distribution charging every point. 
Let $X_0$ be the root of the $\bp$-tree $\TT^{\bp}$ 
and $X_1$ be random with law $\bp$ independent of $\TT^{\bp}$. We
consider $X_1$ as a second root, and call the path $X_0=c_1,c_2,\ldots,c_K=X_1$
from $X_0$ to $X_1$ the {\em spine}. Deleting the edges
$\{c_1,c_2\},\{c_2, c_3\},\ldots$ splits $\TT^{\bp}$ into subtrees rooted at
$c_1,c_2,\ldots, c_K$, which we call $\TT_{c_1},\ldots,\TT_{c_K}$.
Orient the edges of these trees by making them point towards the root. 
Now let $X_2,X_3,\ldots$ be an i.i.d. $\bq$-sample independent of $\TT^{\bp}$. 
Consider the following procedure. 
\begin{itemize}
\item
Let
$\tau_1=2$ and $k_1$ be such that $\TT_{c_{k_1}}$ contains $X_2$.  Bind the
trees $\TT_{c_1},\ldots, \TT_{c_{k_1}}$ by adding oriented edges $c_1\to
c_2\to\ldots\to c_{k_1}\to c_1$.  Let $\cyc_1=\{c_1,\ldots,c_{k_1}\}$ and
$\com_1=\cup_{1\leq i\leq k_1}\TT_{c_i}$.
\item
Given $\tau_i,k_i,\cyc_i,\com_i,1\leq i\leq j$ as long as $\cup_{1\leq i\leq
j}\com_i\neq[n]$, let $\tau_{j+1}=\inf\{k:X_k\notin\cup_{1\leq i\leq
j}\com_i\}$ and $k_{j+1}$ be such that $\TT_{c_{k_{j+1}}}$ contains
$X_{\tau_{j+1}}$.  Then add edges $c_{k_j+1}\to c_{k_j+2}\to\ldots\to
c_{k_{j+1}}\to c_{k_j+1}$, let $\cyc_{j+1}=\{c_{k_j+1},\ldots,c_{k_{j+1}}\}$
, $\com_{j+1}=\cup_{k_j+1\leq i\leq k_{j+1}}\TT_{c_i}$. 
\end{itemize}
When it terminates, say at stage  $r$, the procedure yields a digraph with $r$
connected  components  $\com_1,\ldots,\com_r$,  and  each  component  contains
exactly one cycle of the form $c_{k_j+1}\to \ldots\to c_{k_{j+1}}\to
c_{k_j+1}$. Let $J(\TT^{\bp},X_i,i\geq 1)$ 
be the mapping whose digraph equals the one
given by the procedure. Then, as an easy variation of 
\cite[Proposition 1]{almpimap02}, 

\begin{prp}\label{prpjoyal}
The random mapping $J(\TT^{\bp},X_i,i\geq 1)$ 
is a $\bp$-mapping, and the order on its
basins $\com_1,\com_2,\ldots,\com_r$ induced by the algorithm is $\bq$-biased
order. 
\end{prp}

\subsection{Consequences for associated walks}

From now on, let $\TT^{\bp}$ be a $\bp$-tree, and $H^{\bp}_{\bw}$ be its associated
height process. Let $v_1,v_2,\ldots,v_n$ be the vertices of $\TT^{\bp}$ in depth-first
order, and let $S_{\bw}$ be the linear interpolation between points 
$((\sum_{1\leq j\leq i}w_j,i/n),0\leq i\leq n)$. Given a random variable 
$U$ uniform on $[0,1]$ and independent of $H^{\bp}$, let $X_1=X_1(U)$ be the
vertex  that  is  visited  by  the  walk  at  time  $U^{\bw}=S_{\bw}^{-1}\circ
S_{\bp}(U)$, so this vertex is a $\bp$-distributed random variable independent 
of  $\TT^{\bp}$. We also  let $X_2,X_3,\ldots$  be 
an independent $\bq$-sample, independent of $\TT^{\bp},U$. 
Let $M=J(\TT^{\bp},X_i,i\geq 1)$ be the
$\bp$-mapping associated to $\TT^{\bp}$ by the Joyal correspondence. We will
prove Theorem \ref{T1}(i) by showing that the mapping-walk associated to $M$
converges in law to $Z^{\sth}$. 

Consider the slight variation of the process $\un{H}^{\bp}_{\bw}(u)$:
$$
K^{\bp}_{\bw}(u)(s)=
\left\{\begin{array}{ll}
\un{H}^{\bp}_{\bw}(u)(s)& \mbox{ if } s \mbox{ is not a 
time when a vertex of the spine is visited}\\
\un{H}^{\bp}_{\bw}(u)(s)+1& \mbox{ else.}
\end{array}\right. 
$$
This process thus ``lifts'' the heights of the vertices of the spine by
$1$. Recall from the proof of 
\cite[Lemma 3]{almpimap02} (with a slightly more general context that 
incorporates the weights $\bw$) that these vertices are visited precisely
at the times for
which the reversed pre-minimum process $s\mapsto\un{H}^{\bp}_{(u-s)-}(u)$
jumps downward, so in $K^{\bp}_{\bw}(u)$ we just delay these jumps by the
corresponding $\bw$-mass of the vertex. What we now call ``excursion'' or 
generalized excursion of $H^{\bp}_{\bw}$ above $K^{\bp}_{\bw}(u)$ is just 
the same as before, that is a recentered portion of the path of 
$H^{\bp}_{\bw}$ on a flat interval of $K^{\bp}_{\bw}(u)$, with the convention 
that two excursions on two flat intervals with same heights (here and below the term 
``height'' refers to the flat intervals) are 
merged together as a single generalized excursion. By contrast with the above, 
these excursions may take negative values, but only at times when cyclic 
vertices are visited, where the excursions' value is $-1$. As above, let
$\wt{\bJ}^{u}(H^{\bp}_{\bw})$ be the process obtained by merging the 
excursions of $H^{\bp}_{\bw}$ above $K^{\bp}_{\bw}(u)$ in increasing 
order of height. A slight variation of \cite[Lemma 3]{almpimap02} gives

\begin{lmm}\label{lemj}
$$\wt{\bJ}^{U^{\bw}}(H^{\bp}_{\bw})=H^M_{\bw}-1.$$
\end{lmm}
Notice in  particular that $H^M_{\bw}$  is a functional  of $\TT^{\bp}$ and $X_1(U)$
alone, and does not depend on $X_2,X_3,\ldots$. 

\noindent{\bf Proof of Theorem \ref{T1}. }
Let $\bp_n$ satisfy (\ref{regime})
with finite-length limit $\btheta$. 
We use Theorem \ref{thptree} and
Skorokhod's representation theorem, so we suppose that the convergence of
$\spn H^{\bp_n}\to H^{\sth}$ (either in $*$-topology or Skorokhod topology
according to the hypotheses) is almost-sure, as well as the convergence of
$S_{\bp_n},S_{\bw_n},S_{\bq_n}$ to the identity. We also suppose that the
convergence of Theorem \ref{thptree} (ii) is almost-sure. 
%The variable $U$ is taken to be the
%same for every $n$  up to enlarging the probability space. 

Fix $\epsilon>0$. 
For (Lebesgue) almost-every $u\in[0,1]$, $u$ is not a local minimum of
$H^{\sth}$ on the right or on the left.  Fix such a $u$.  Since
$u^{\bw_n}:=S^{-1}_{\bw_n}\circ S_{\bp_n}(u)\to u$ as $n\to\infty$, it is
easily checked that for any $\eta>0$ and $n>N_1$ large enough, the processes
$\un{H}^{\sth}(u^{\bw_n})$ and $\un{H}^{\sth}(u)$ (resp.\
$K^{\bp_n}_{\bw_n}(u^{\bw_n})$ and $K^{\bp_n}_{\bw_n}(u)$) co\"\i ncide
outside  the  interval $(u-\eta,u+\eta)$.  Let  $\eps_1,\eps_2,\ldots$ be  the
generalized excursions of $H^{\sth}$ above $\un{H}^{\sth}(u)$, ranked by
decreasing  order of their  durations $l_1,l_2,\ldots$,  call $h_1,h_2,\ldots$
the corresponding (pairwise distinct) heights. Let $\alpha>0$ be
such that
$\omega(h):=\sup_{h\in[-\alpha,\alpha]}||H^{\sth}_{\cdot+h}-
H^{\sth}_{\cdot}||_{\infty}<\epsilon/3$. Notice that 
for $n>N_2$ large enough, we also have 
$\omega_n(h):=\spn\sup_{h\in[-\alpha,\alpha]}||H^{\bp_n}_{\bw_n}(\cdot+h)-
H^{\bp_n}_{\bw_n}(\cdot)||_{\infty}\leq \eps/2$. 
Next, take $k$ such that $\sum_{i=1}^kl_i\geq 1-\alpha/2$, and
choose $\eta<\alpha/4$ such that
none of the intervals of constancy of $\un{H}^{\sth}(u)$ corresponding to
these $k$ excursions intersect $(u-\eta,u+\eta)$. 

Next, consider hypothesis (i) of 
Theorem \ref{T1}.  If $[a,b]$ is an interval of constancy of 
$\un{H}^{\sth}(u^{\bw_n})$ (or $\un{H}^{\sth}(u)$) 
not intersecting $(u-\eta,u+\eta)$,
then there exists for $n$ large enough a constancy interval
of $K^{\bp_n}_{\bw_n}(u)$, which we denote by $[a^n,b^n]$, such that 
$(a^n,b^n)\to (a,b)$, implying by Theorem \ref{thptree}(i) that 
$$(\spn (H^{\bp_n}_{\bw_n}(a^n+s)-H^{\bp_n}_{\bw_n}(a^n)),0\leq s\leq b^n-a^n)
\to (H^{\sth}(a+s)-H^{\sth}(a),0\leq s\leq b-a)$$
uniformly.  Moreover, for $u$ as chosen above, if $u\in(t_i,T_i)$ (notice
$u=T_i$ or $u=t_i$ is not possible) then there
exists some $t_{ij},T_{ij}$ with $t_{ij}<u<T_{ij}$. 
Thus, for  such $u$  and as  a consequence of  Theorem \ref{thptree}  (ii), if
there  exists a  second such  flat interval  $[c,d]$ with  same height  as the
initial one (with say $b<c$), then there
also exists a constancy interval $[c^n,d^n]$ of $K^{\bp_n}_{\bw_n}(u)$ with 
$(c^n,d^n)\to(c,d)$, with the {\bf same} height as the first one. 
Therefore, these two intervals do merge to form the interval of a 
generalized excursion of 
$\spn H^{\bp_n}_{\bw_n}$ above $\spn K^{\bp_n}_{\bw_n}(u)$ with length 
$(b^n-a^n)+(d^n-c^n)$, that converges uniformly to the generalized 
excursion of $H^{\sth}$ above $\un{H}^{\sth}(u)$ with height $H^{\sth}_{a}$ 
and duration $(b-a)+(d-c)$.  As a conclusion, one has $\eps_i^n\to \eps_i$
uniformly for every $1\leq i\leq k$, where $\eps_i^n$ is the generalized 
excursion of $\spn H^{\bp_n}_{\bw_n}$ above $\spn 
K^{\bp_n}_{\bw_n}(u^{\bw_n})$ with $i$-th largest
duration $l_i^n$. Call $h_i^n$ its height. 

Now $(h^n_1,\ldots,h^n_k)\to (h_1,\ldots,h_k)$, 
and $\sum_{1\leq i\leq k}|l^n_i-l_i| \to 0$ as $n\to\infty$. 
%Fix $\epsilon>0$. Choose $\alpha>0$ such that the uniform norm of 
%$H^{\sth}_{\cdot+h}-H^{\sth}_{\cdot}$ is $\leq\epsilon/2$ for every 
%$h\in[-\alpha,\alpha]$, and choose $k$ such that both $\sum_{1\leq i\leq
%k}l_i>1-\alpha/2$ and $||\eps_i||_{\infty}<\epsilon/4$ for $i>k$ (this is possible
%by continuity of $H^{\sth}$). 
Thus, if $n>\max(N_1,N_2)$ is also chosen so that 
\begin{itemize}
%\item $||H^{\bp_n}_{\bw_n}-H^{\sth}||_{\infty}<\epsilon/4$ and
%$||K^{\bp_n}_{\bw_n}(u)-\un{H}^{\sth}(u)||_{\infty}<\epsilon/4$,
\item $\sum_{1\leq i\leq k}|l^n_i-l_i| \leq \alpha/2$,
\item $h^n_1,\ldots,h^n_k$  are in the same order  as $h_1,\ldots,h_k$ (recall
these are pairwise distinct),
\item $\sup_{1\leq i\leq k}||\eps^n_i-\eps_i||_{\infty}<\epsilon/2$, 
\end{itemize}
then necessarily, the uniform distance between 
$\spn \wt{\bJ}^{u^{\bw_n}}(H^{\bp_n}_{\bw_n})$ and 
$\bJ^{u^{\bw_n}}(H^{\sth})$ is at most $\epsilon$. Indeed, for $x\in [0,1]$, 
if $x\in(g_i^n,d_i^n)\cap(g_i,d_i)$
for some $i\leq k$, then
$$|\spn\wt{\bJ}^{u^{\bw_n}}(H^{\bp_n}_{\bw_n})(x)-
\bJ^{u^{\bw_n}}(H^{\sth})(x)| \leq
||\eps_i-\eps_i^n||_{\infty}+
\sup_{|h|<\alpha}||\eps_i(\cdot)-\eps_i(\cdot+h)||_{\infty}\leq \epsilon,$$
and else the value taken by this difference does not exceed 
$\omega(h)+\omega_n(h)\leq  \epsilon$ because  there must  be a  zero  of both
processes at distance $< \alpha$ from $x$. 
Apply this to
$u=U$, which a.s.\ does not belong to  the set of local minima (on the left or
on the right) of $H^{\sth}$. 
Using Lemma \ref{lemj} establishes the assertion of (i).

The case (iv) of $*$-convergence follows the same lines as in 
\cite[Lemma 2]{almpimap02}. We suppose  up to extracting subsequences that $\spn H^{\bp_n}$
can be written as $g_n+h_n$ with $g_n$ converging uniformly to $H^{\sth}$ and 
$h_n(u)=0$ ultimately for almost-every $u$. 
Then, up  to modifying  slightly the constancy  intervals of  $K^{\bp_n}$, the
same result as above holds for $g_n$, 
so this proves that $\spn \bJ^U(H^{\bp_n})$
converges to $Z^{\sth}$ in probability for the $*$-metric. 

Points (ii,iii) in Theorem \ref{T1} then follow the same lines as in the
proof of \cite[Theorem 1]{almpimap02}. We give some details for (ii). Let
$U_2,U_3,\ldots$ be uniform independent random variables, independent of 
$H^{\sth},U,(H^{\bp_n},n\geq 1)$.  Let $U_i^{\bq_n}=S_{\bw_n}^{-1}\circ
S_{\bq_n}(U_i)$ for $i\geq2$.  Recall that the walk $H^{M_n}_{\bw_n}$ 
can  be defined  using only $H^{\bp_n},U$,  so we  are allowed  to  make the
following choice for $X_2,X_3,\ldots$: we let $X_i$ be the vertex
visited by $H^{M_n}_{\bw_n}$ at time $U_i^{\bq}$.  Therefore, the marks
$D^{M_n}_{\bw_n}(i)$ are obtained recursively as follows: let $v$ be the
vertex visited by the first $U^{\bq}_j>D^{M_n}_{\bw_n}(i)$, then 
$D^{M_n}_{\bw_n}(i+1)$ is the first time when a cyclic point is visited
strictly after $v$, i.e.\ at the right end of the generalized
excursion of $H^{M_n}_{\bw_n}$ straddling this $U^{\bq}_j$. 
Passing to the limit, we find that $(D^{M_n}_{\bw_n}(i),1\leq i\leq j)$ 
converges a.s.\ to $(D'_i,1\leq i\leq j)$ defined
recursively by: $D'_{i+1}$ is the first point of $\{d_1,d_2,d_3,\ldots \}$ that occurs
after the  first $U_j>D'_i$.  It is easy  to see that this  defines a sequence
with the same law as $D_i,i\geq 1$. \cq

\section{Inhomogeneous continuum random tree interpretation}\label{sec-ICRT}

Let us briefly introduce the details of the 
limiting ICRT's {\em stick-breaking} construction
\cite{jpmc97b,jpda97ebac}.  Let $\btheta = (\theta_0, \theta_1, \theta_2, \ldots )$ satisfy $\sum_{i \geq 0} \theta_i^2 = 1$.  Consider a Poisson process $(U_j,V_j),j\geq 1$ on
the first octant $\mathbb{O}=\{(x,y)\in\R^2:0\leq y\leq x\}$, with intensity
$\theta_0^2$ per unit area. For each $i\geq 1$ consider also homogeneous 
Poisson processes $(\xi_{i,j},j\geq 1)$ with intensity $\theta_i$ per unit
length,
and suppose these processes are
independent,  and  independent  of  the  first  Poisson  process.  The  points
of $\R_+$ that are either equal to some $U_i,i\geq 1$ or some $\xi_{i,j},j\geq
2$ will be called {\em cutpoints}. To a cutpoint $\eta$ we associate a {\em
joinpoint} $\eta^*$: if $\eta$ is of the form $U_i$, let $\eta^*=V_i$, while
if $\eta=\xi_{i,j}$ for some $i\geq 1,j\geq 2$, we let
$\eta^*=\xi_{i,1}$. Under the hypothesis $\sum_i\theta_i^2$, one shows that we
may order the cutpoints as $0<\eta_1<\eta_2,\ldots$ almost-surely.  We build
recursively a consistent family of trees whose edges are line-segments by first
letting $\TT^{\sth}_1$ be the segment $[0,\eta_1]$ rooted at $0$, and then,
given $\TT^{\sth}_J$, by attaching the left-end of the segment
$(\eta_J,\eta_{J+1}]$ at the corresponding joinpoint $\eta_J^*$, which has
been already placed somewhere on $\TT^{\sth}_J$.  Further, we relabel the
joinpoints of the form $\xi_{i,1}$ as $i$, and we relabel the leaves
$\eta_1,\eta_2,\ldots$ as $1+,2+,\ldots$. When all the branches are attached,
we obtain a random metric space whose completion we call $\TT^{\sth}$ (it can
therefore be interpreted as the completion of a special metrization of
$[0,\infty)$). We  let $[[v,w]]$ be the  only injective path from  $v$ to $w$,
and $]]v,w]]=[[v,w]]\setminus\{v\}$. 

Together with the ICRT comes one natural measure, which is 
the length measure inherited from Lebesgue measure 
on $[0,\infty)$. When $\btheta$ satisfies the further hypothesis $\theta_0>0$
or $\sum_i\theta_i=\infty$, the tree can be endowed with another measure
$\mu$, which is a probability measure obtained as the weak limit of the
empirical distribution $\mu_J$ on the leaves $1+,2+,\ldots,J+$ as 
$J\to\infty$. We call $\mu$ the {\em mass measure}. 

If $\btheta\in\Thfinite$, it
has been shown in \cite{almpiep} that $H^{\sth}$ is the 
{\em exploration process} of $\TT^{\sth}$. 
To explain what this means, note first that 
$H^{\sth}$
 induces  a  special
pseudo-metric on $[0,1]$ by letting
$$d(u,v)=H^{\sth}_u+H^{\sth}_v-2\inf_{w\in[u,v]}H^{\sth}_w.$$
It turns out that the quotient space $\TT$ obtained by
identifying points of $[0,1]$ at distance $0$ has the same ``law'' as
$\TT^{\sth}$, where the mass measure is the measure on the quotient induced by
Lebesgue measure on $[0,1]$. %   Precisely, the result states that if
Precisely,
\begin{thm}[\cite{almpiep}]
\label{T9}
If $U_1,\ldots,U_J$ are independent uniform variables on $[0,1]$, independent of
$H^{\sth}$, then the subtree of $\TT^{\sth}$ spanned by the (equivalence classes of
the) $U_i$'s has the same law as $\TT^{\sth}_J$. \end{thm}    % process $H^{\sth}$
Conceptually, the stick-breaking construction provides an
``algorithmic construction" of the ICRT, whereas the process $H^{\sth}$
plays a r\^ole similar to that of Brownian excursion
in our methodology described in point (ii) in the introduction.

We now show how some consequences of our main theorem can be formulated in
terms of the stick-breaking construction of the ICRT.  % as it might be interesting
 For $v\in\TT^{\sth}$, let ${\tt junc}(v)$ be the
branchpoint between $v$ and $1+$. 
Define recursively a sequence
$0=c_0,c_1,\ldots$ of vertices of the spine $[[{\rm root},1+]]$ with increasing
heights recursively using the rule
\begin{quote}
Given $c_j$ let $k_{j+1}+$ be the first leaf of $\{2+,3+,4+,\ldots \}$ with ${\tt
junc}(k_{j+1}+)\notin[[{\rm root},c_j]]$ and let $c_{j+1}={\tt junc}(k_{j+1}+)$. 
\end{quote}

\begin{crl}\label{crlicrt}
Under regime (\ref{regime}) with limiting $\btheta\in\Thfinite$,
\begin{eqnarray*}
\lefteqn{(\bp_n(\com_j(M_n)),\spn\card(\cyc_j(M_n)),j\geq 1)}\\
&\to& 
\left(\lim_{k\to \infty}\frac{1}{k}\card\{1\leq i\leq k:{\tt
junc}(i+)\in]]c_{j-1},c_j]]\}, \HT(c_j)-\HT(c_{j-1}),j\geq 1\right)
\end{eqnarray*}
\end{crl}

\proof
The $n \to \infty$ limit of the left side is (by Corollary \ref{C1})
the law of
\begin{equation}
(D_j-D_{j-1},L^{\sth}_{D_j}-L^{\sth}_{D_{j-1}},j\geq 1) 
 . \label{newD} \end{equation}
By the description of $\mu$ as the $k \to \infty$ limit
of the empirical distribution on leaves $\{1+,2+,\ldots,k+\}$,
the $k \to \infty$ limit of the right side % (by Corollary \ref{C1})
of Corollary \ref{crlicrt} becomes
\begin{equation}
(\mu\{v\in \TT^{\sth}:{\tt
junc}(v)\in]]c_{j-1},c_j]]\},\HT(c_j)-\HT(c_{j-1}),j\geq  1)  % We  conclude by
 . \label{newmu} \end{equation}
So the issue is to show equality in law of (\ref{newD})
and (\ref{newmu}).
But Theorem \ref{T9} identifies the law (\ref{newmu}) with the law
\begin{equation}
(\mathrm{Leb}\{v\in (0,1):{\tt
junc}(v)\in]]c_{j-1},c_j]]\},H^{\sth}_{c_j}-H^{\sth}_{c_{j-1}},j\geq  1)  % We  conclude by
  \label{newH} \end{equation}
where the quantities involved can be redefined as follows.
Take $U_1,U_2,U_3,\ldots $ uniform on $(0,1)$, independent of $H^{\sth}$.
Let ${\tt junc}(v)$ be the point at which 
$\inf_{[v,U_1]} H^{\sth}_{\cdot}$ or
$\inf_{[U_1,v]} H^{\sth}_{\cdot}$ is attained.
Given $c_j$, let $c_{j+1} = {\tt junc}(U^\prime)$
where $U^\prime$ is the first of 
$\{U_2,U_3,U_4,\ldots\}$
such that
$H^{\sth}_{{\tt junc}(U^\prime)} > H^{\sth}_{c_j}$.

On the other hand,
 $D_1$ is by definition equal in law to the sum of the lengths
of the  generalized excursions  of $H^{\sth}$ above  $\un{H}^{\sth}(U_1)$ 
whose
heights are less than or equal to that of the excursion containing an
independent uniform $U_2$, while $L^{\sth}_{D_1}$ is the height of the 
corresponding excursion. Recursively, $D_{j+1}-D_j$ is equal in law to the sum of the
durations of the excursions with heights between the height of the previously
explored excursions (strictly) and the height of the excursion straddling the
first $U_i$ that falls in an excursion interval with height larger than the
previous ones; $L^{\sth}_{D_j}-L^{\sth}_{D_{j-1}}$ is then the difference of these heights. 
This identifies the law (\ref{newD}) with the law (\ref{newH}). \cq

%Assimilating  the $U_i$'s  to  $\mu$-distributed leaves  of $\TT^{\sth}$,  and
%further, assimilating these leaves to $1+,2+,\ldots$ ({\tt yyy here I might be
%missing some argument}), 
%we obtain that
%$(D_j-D_{j-1},L^{\sth}_{D_j}-L^{\sth}_{D_{j-1}},j\geq 1)$ 
%has same law as 
%$(\mu\{v\in \TT^{\sth}:{\rm
%junc}(v)\in]]c_{j-1},c_j]]\},\HT(c_j)-\HT(c_{j-1},j\geq  1)$. We  conclude by
%the strong law of large numbers. \cq

\rem 
Corollary \ref{crlicrt} could alternatively be proved, for more general limit regimes, by an argument based directly on the Joyal
correspondence, without using the detour through exploration processes. % and for more
%general regimes, namely for $\btheta\notin\Thfinite$. 
%$\bullet$ Cases where $\btheta$ satifies $\theta_0=0,\sum_i\theta_i<\infty$
%are less interesting as the corresponding ICRT's do not have a fine
%structure. For example,  if $\theta_1=1$ then $\TT^{\sth}$ is  a collection of
%iid. exponential($1$)-long sticks connected to  the same point by one of their
%ends.  This corresponds to a limiting  mapping with only two basins, the first
%one being a pure cycle  with exponential($1$) length, and the other containing
%a  cycle with  exponential($1$)  length  in which  exactly  one cyclic  points
%attracts the  trajectories of all non-cyclic points  (these trajectories being
%unconnected lines with exponential($1$) lengths). 

\section{Final remarks}
The regimes (\ref{regime}) are basically the only
possible ones, if we require a limit distribution for the number $|{\cal C}(M_n)|$
of cyclic vertices.
\begin{lmm}
If $c_n(|{\cal C}(M_n)|-d_n)$ converges in law to some non-trivial
distribution on $\R_+$ for some renormalizing sequences $c,d$, then there
exists $\btheta$ such that $\bp$ satisfies (\ref{regime}) up to
elementary rescaling, that is, there exists $\alpha\in(0,\infty)$
and $\beta\in\R$ such that $c_n/\spn\to\alpha$ and $c_nd_n\to
\beta$.
\end{lmm} 

This lemma is a direct consequence of
\cite[Theorem  4]{jpmc97b}  and  of  Proposition \ref{prpjoyal},  which
implies that 
the number of cyclic points of a $\bp$-mapping has same distribution as one
plus the distance from the root to a $\bp$-sampled vertex of a $\bp$-tree.

\end{document}